\documentclass[a4paper,12pt]{amsart}
\usepackage{amsfonts}
\usepackage{amssymb}
\usepackage{ifthen}
\usepackage{graphicx}
\usepackage{color}
\nonstopmode \numberwithin{equation}{section}
\newtheorem{thm}{Theorem}%[section]
\newtheorem{lem}{Lemma}%[section]
\newtheorem{cor}{Corollary}%[section]
\newtheorem{prop}{Proposition}%[section]

\newtheorem{conj}{Conjecture}

\theoremstyle{definition}
\newtheorem{defn}{Definition}%[section]
\newtheorem{example}{Example}%[section]

\newtheorem{ques}{Question}
\newtheorem{rem}{Remark}

\newtheorem{rems}{Remarks}

\newcounter {own}
\def\theown {\thesection  .\arabic{own}}

\newenvironment{pf}[1][]{%
 \vskip 3mm
 \noindent
 \ifthenelse{\equal{#1}{}}%
  {{\slshape Proof. }}%
  {{\slshape #1.} }%
 }%
{\qed\bigskip}

\newcounter{alphabet}
\newcounter{tmp}
\newenvironment{Thm}[1][]{\refstepcounter{alphabet}%
\medskip%
\noindent%
{\bf Theorem \Alph{alphabet}}%
\ifthenelse{\equal{#1}{}}{}{ (#1)}%
{\bf .} \itshape}{\vskip 8pt}
\makeatletter
\newcommand{\Ref}[1]{\@ifundefined{r@#1}{}{\setcounter{tmp}{\ref{#1}}\Alph{tmp}}}
\makeatother
\newenvironment{Lem}[1][]{\refstepcounter{alphabet}%
\medskip%
\noindent%
{\bf Lemma \Alph{alphabet}}%
{\bf .} \itshape}{\vskip 8pt}

\newcommand{\N}{{\mathbb N}}

\newcommand{\T}{{\mathbb{T}}}

\newcommand{\Aut}{{\operatorname{Aut}}}

\def\be{\begin{equation}}
\def\ee{\end{equation}}

\newcommand{\bee}{\begin{enumerate}}
\newcommand{\eee}{\end{enumerate}}

\newcommand{\blem}{\begin{lem}}
\newcommand{\elem}{\end{lem}}
\newcommand{\bthm}{\begin{thm}}
\newcommand{\ethm}{\end{thm}}
\newcommand{\bcor}{\begin{cor}}
\newcommand{\ecor}{\end{cor}}
\newcommand{\beg}{\begin{example}}
\newcommand{\eeg}{\end{example}}
\newcommand{\begs}{\begin{examples}}
\newcommand{\eegs}{\end{examples}}
\newcommand{\bdefn}{\begin{defn}}
\newcommand{\edefn}{\end{defn}}
\newcommand{\bprob}{\begin{prob}}
\newcommand{\eprob}{\end{prob}}
\newcommand{\bei}{\begin{itemize}}
\newcommand{\eei}{\end{itemize}}
\newcommand{\bqn}{\begin{ques}}
\newcommand{\eqn}{\end{ques}}
\newcommand{\bcon}{\begin{conj}}
\newcommand{\econ}{\end{conj}}
\newcommand{\bcons}{\begin{conjs}}
\newcommand{\econs}{\end{conjs}}
\newcommand{\bprop}{\begin{prop}}
\newcommand{\eprop}{\end{prop}}
\newcommand{\brem}{\begin{rem}}
\newcommand{\erem}{\end{rem}}
\newcommand{\brems}{\begin{rems}}
\newcommand{\erems}{\end{rems}}
\newcommand{\bo}{\begin{obser}}
\newcommand{\eo}{\end{obser}}
\newcommand{\bos}{\begin{obsers}}
\newcommand{\eos}{\end{obsers}}
\newcommand{\bpf}{\begin{pf}}
\newcommand{\epf}{\end{pf}}
\newcommand{\ba}{\begin{array}}
\newcommand{\ea}{\end{array}}
\newcommand{\beq}{\begin{eqnarray}}
\newcommand{\beqq}{\begin{eqnarray*}}
\newcommand{\eeq}{\end{eqnarray}}
\newcommand{\eeqq}{\end{eqnarray*}}

\newcommand{\Ra}{\Rightarrow}

\newcommand{\ra}{\rightarrow}

\newcommand{\ds}{\displaystyle}

\newcounter{minutes}
\divide\time by 60
\newcounter{hours}
\multiply\time by 60 \addtocounter{minutes}{-\time}
\begin{document}
\bibliographystyle{amsplain}
\title[Composition operators on Hardy spaces of the homogenous rooted trees]
{Composition operators on Hardy spaces of the homogenous rooted trees}

%[Composition operators on the discrete Hardy space on homogenous trees-II]
%{Composition operators on the discrete Hardy space on homogenous trees-II}

%%%=========================================================================
%%\thanks{%$^\dagger$
%File:~\jobname .tex,
%          printed: \number\day-\number\month-\number\year,
%          \thehours.\ifnum\theminutes<10{0}\fi\theminutes}
%%%=========================================================================
%
%\author{S. Ponnusamy $^\dagger $}

\author{Perumal Muthukumar}
\address{P. Muthukumar, Indian Statistical Institute,
Statistics and Mathematics Unit, 8th Mile, Mysore Road,
Bangalore, 560059, India.}
\email{pmuthumaths@gmail.com}

\author{Saminathan Ponnusamy}
\address{S. Ponnusamy, Department of Mathematics,
Indian Institute of Technology Madras, Chennai, 600 036, India.}
\email{samy@iitm.ac.in}


%\date{\today

\begin{abstract}
In \cite{CO-Tp-spaces}, the present authors initiated the study of composition operators on discrete analogue of generalized
Hardy space $\mathbb{T}_{p}$ defined on a homogeneous rooted tree. In this article, we give equivalent conditions for
the composition operator $C_\phi$ to be bounded on $\mathbb{T}_{p}$ and on $\mathbb{T}_{p,0}$
spaces and compute their operator norm. We also characterize
invertible composition operators as well as isometric composition operators on $\mathbb{T}_{p}$ and on $\mathbb{T}_{p,0}$
spaces. Also, we discuss the compactness of
$C_\phi$ on $\mathbb{T}_{p}$ and finally prove there are no compact composition operators
on $\mathbb{T}_{p,0}$ spaces.
\end{abstract}
\thanks{
%%=========================================================================
File:~\jobname .tex,
          printed: \number\day-\number\month-\number\year,
          \thehours.\ifnum\theminutes<10{0}\fi\theminutes
%%=========================================================================
}
\maketitle
\pagestyle{myheadings}
\markboth{P. Muthukumar and S. Ponnusamy}{Composition Operators on the Discrete Hardy Space on Homogenous Trees}

\section{Introduction}\label{MP5Sec1}

Let $X$ be a Banach space of complex valued functions  on a nonempty set $\Omega$.
The composition operator $C_\phi$ induced by a self-map $\phi$ of $\Omega$ is defined as
$$ C_\phi(f)=f\circ\phi ~\mbox{ for all  $f\in X$}.
$$
The study of composition operators on analytic function spaces has a rich history. The typical choices
for $X$ are spaces of analytic functions on the unit disk
such as the Hardy spaces, the Bergman spaces, the Bloch space, or the Dirichlet spaces. We refer to the book of Cowen and MacCluer \cite{Cowen:Book}
for composition operators defined on various spaces of analytic functions on the
unit disk, whereas the book of Shapiro \cite{Shapiro:Book} is devoted mainly to composition operators on Hardy spaces.
The composition operators on various measure spaces are discussed in the book of Singh and Manhas \cite{Rksingh:Book}.
These books bring together many well-developed aspects of the subject along with several open problems.
The systematic study of operator theory on discrete structure specially on infinite trees has been the subject of several recent papers
\cite{Colonna-MO-5,Colonna-MO-3,Colonna-MO-2,Colonna-CO,Colonna-MO-6,Allen-MO-7,Allen-CO-2,Colonna-MO-1,Colonna-MO-4,Colonna-Toeplitz,MP-Tp-spaces,CO-Tp-spaces}.

Discrete function spaces are mostly defined to be analogs of analytic function spaces (cf. \cite{CohnColo-96}).
Multiplication and composition operators are mainly considered on discrete function spaces. The basic questions
such as boundedness, compactness, estimates for operator norm and essential norm, isometry and spectrum were
considered for multiplication operators between various discrete function spaces on infinite tree such as
Lipschitz space, weighted Lipschitz space and iterated logarithmic Lipschitz spaces.
See  \cite{Colonna-MO-5,Colonna-MO-3,Colonna-MO-2,Colonna-MO-6,Colonna-MO-1,Colonna-MO-4} for more details.

The study of composition operators on discrete function space was first initiated by Colonna et al. \cite{Colonna-CO}.
In that paper  the Lipschitz space of a tree was investigated. Multiplication and composition operators on weighted
Banach spaces of an infinite tree were considered in
\cite{Allen-MO-7,Allen-CO-2}, respectively. Recently, some classes of operators including Toeplitz operators
with symbol from the Lipschitz space of a tree were considered in \cite{Colonna-Toeplitz}.
In \cite{MP-Tp-spaces}, the present authors defined discrete analogue ($\mathbb{T}_{p}$) of generalized Hardy spaces on
homogeneous rooted tree and studied multiplication operators on them. Study of
composition operators on $\mathbb{T}_{p}$ spaces were initiated by the present authors in \cite{CO-Tp-spaces}.

In this article,  we continue the study of composition operators on $\mathbb{T}_{p}$ spaces.
We refer to  Section \ref{MP5Sec2} for preliminaries about $\mathbb{T}_{p}$ and $\mathbb{T}_{p,0}$ spaces.
In Subsection \ref{MP5Subsec1}, we give equivalent conditions for
the composition operator $C_\phi$ to be bounded on various $\mathbb{T}_{p}$ spaces and compute their operator norms.
Subsection \ref{MP5Subsec2} is devoted to the study of composition operators induced by special symbols such as
injective and multivalent maps. In Subsection \ref{MP5Subsec3}, we  discuss bounded composition operators on
$\mathbb{T}_{p,0}$ and their norm estimates. In Sections \ref{MP5Sec4} and \ref{MP5Sec5},
we characterize invertible composition operators and isometric
composition operators on various $\mathbb{T}_{p}$ and $\mathbb{T}_{p,0}$ spaces.
Finally, in Section \ref{MP5Sec6},  we present some results about compactness of $C_\phi$ on $\mathbb{T}_{p}$ and
 prove that there are no compact composition operators
on $\mathbb{T}_{p,0}$ spaces.

\section{Preliminaries and Lemmas}\label{MP5Sec2}
To make the paper self-contained we recall some basic definitions. More details can be found in
standard texts on this subject (cf. \cite{Book:graph}).

Let $G=(V,E)$ be a graph such that $E\subseteq V \times V$, where the elements of the sets $V$ and $E$ are called vertices
and edges of the graph $G$, respectively. We shall not always distinguish between a graph and its vertex set and so,
we may write $x \in G$ (rather than $x \in V$) and by a function defined on a graph, we mean a function defined on its vertices.
Similarly, a self-map of a graph is a function defined on its vertices to itself.
Two vertices $x,y \in G$  are said to be \textit{neighbours} (denoted by $x\sim y$) if $(x,y)\in E$.

A graph  is said to be \textit{$k$-homogeneous} if every vertices of the graph have exactly $k$ neighbours.
A \textit{finite path} is a nonempty subgraph $P=(V,E)$ of the form $V=\{x_0,x_1,\ldots ,x_k\}$ and
$E=\{(x_0,x_1),(x_1,x_2),\ldots,(x_{k-1},x_k)\}$, where $x_i$'s are distinct. In this case, we call $P$ a path
between $x_0$ and $x_k$. If $P$ is a path between $x_0$ and $x_k ~(k\geq 2)$, then $P$ with an additional
edge $(x_{n},x_{0})$ is called a \textit{cycle}. A nonempty graph is said to be \textit{connected} if
there is a path between any two of its vertices. A connected  graph without cycles is called a \textit{tree}. Thus, any two vertices of a
tree are linked by a unique path. The \textit{distance} between any two vertex of a
tree is the number of edges in the unique path connecting them. Sometimes it is convenient to consider one vertex of a tree as
special; such a vertex is then called the root of this tree.

A tree with fixed root $\textsl{o}$ is called a \textit{rooted tree}.
If $T$ is a rooted tree with root $\textsl{o}$, then $|v|$ denotes the distance between the root $\textsl{o}$ and the vertex $v$.
Further, the \textit{parent} (denoted by $v^-$) of a vertex $v$, which is not a root, is the unique vertex $w\in T$ such that
$w\sim v$ and $|w|=|v|-1$. In this case, $v$ is called \textit{child} of $w$.

Throughout the paper, unless otherwise stated explicitly, $T$ denotes a homogeneous rooted tree (hence an infinite graph),
$\phi$ denotes a self-map of $T$, $\mathbb{N}=\{1,2,\ldots \}$  and $\mathbb{N}_0=\mathbb{N}\cup \{0\}$.

%For $p\in(0,\infty]$, the {\it  Hardy space}
%$H^{p}$ consists of all those analytic functions $f:\mathbb{D}\rightarrow\mathbb{C}$ such that
%$\|f\|_{p}<\infty$,
%where
%$$
%\|f\|_{p}=\sup_{0\leq r<1}M_{p}(r,f)
%$$
%and
%$$M_{p}(r,f)=
%\begin{cases}
%\displaystyle \left(\frac{1}{2\pi}\int_{0}^{2\pi}|f(re^{i\theta})|^{p}\,d\theta\right)^\frac{1}{p}
%& \mbox{if } p\in(0,\infty)\\
%\displaystyle\sup_{|z|=r}|f(z)| &\mbox{if } p=\infty.
%\end{cases}
%$$
%The {\it generalized Hardy space} $H^{p}_{g}$ is defined similarly, upon replacing analytic functions by measurable functions.

As in \cite{MP-Tp-spaces}, for a $(q+1)$-homogeneous tree $T$ rooted at $\textsl{o}$, we define
$$\|f\|_{p}:= \sup\limits_{n\in \mathbb{N}_{0}} M_{p}(n,f),
$$
where $M_{p}(0,f):= |f(\textsl{o})|$ and for every $n\in \mathbb{N}$,
$$%\be\label{MP2eq2}
M_{p}(n,f):=
\left\{
\begin{array}{ll}
\ds \left (\frac{1}{(q+1)q^{n-1}}\sum\limits_{|v|=n}|f(v)|^{p} \right )^{\frac{1}{p}} & \mbox{ if }~ p\in(0,\infty), \\
\max\limits_{|v|=n } |f(v)| & \mbox{ if }~ p=\infty.
\end{array}
\right.
$$
The discrete analogue of the generalized Hardy space, denoted by $\mathbb{T}_{p}$, is then defined by
$$\mathbb{T}_{p}:=\{f\colon T \to\mathbb{C} \, \big |\, \|f\|_{p}<\infty\}
$$
for every $p\in(0,\infty]$.
For the sake of simplicity, we ignore $q$ in the notation of this space.
Similarly, the discrete analogue of the generalized little Hardy space, denoted by $\mathbb{T}_{p,0}$, is defined by
$$\mathbb{T}_{p,0}:=\{f\in\mathbb{T}_{p} :\, \lim\limits_{n\rightarrow\infty} M_{p}(n,f)=0 \}
$$
for every $p\in(0,\infty]$.

Let us fix some notation for the rest of the paper.
Let $T$ be a $(q+1)$-homogeneous tree and $\phi$ denote a self-map of $T$. For $n\in \mathbb{N}_{0}$,
let $D_n$ denote the set of all vertices
$v \in T$ with $|v|=n$ and  denote the number of elements in $D_n$ by $c_n$. Thus,
\be\label{cn}
c_n=\left\{\begin{array}{ll}
(q+1)q^{n-1} & \mbox{ if }~ n\in \mathbb{N}, \\
1 & \mbox{ if }~ n=0.
\end{array} \right.
\ee
For $n\in \mathbb{N}_{0}$ and $w\in T$,
let $N_{\phi}(n,w)$ denote the number of pre-images of $w$ for $\phi$ in $|v|=n$. That is,
$N_{\phi}(n,w)$ is the number of elements in $\{\phi^{-1}(w)\}\bigcap D_n$. Finally,
for each $m$ and $n\in \mathbb{N}_{0}$, $N_{m,n}$ denotes the
maximum of $N_{\phi}(n,w)$ over $|w|=m$. It is obvious that $\sum\limits_{m=0}^\infty N_{m,n}\leq c_n$ for each $n$.

Unless otherwise stated, throughout the discussion,  $\|.\|$ denotes $\|.\|_p$  in $\mathbb{T}_{p}$ spaces. The following results proved
by the present authors in \cite{MP-Tp-spaces} are needed elsewhere.
% in our present investigation.

\begin{Thm}\label{thm:banachp}
{\rm (\cite[Theorems 3.1 and 3.5]{MP-Tp-spaces})}
For $1\leq p\leq\infty$, $\|.\|_{p}$ induces a Banach space structure on the spaces $\mathbb{T}_{p}$ and
$\mathbb{T}_{p,0}$.
\end{Thm}

In \cite{MP-Tp-spaces}, authors raised the question whether $\mathbb{T}_2$ is a Hilbert space or not. The answer is indeed No!.
For example,
Choose two vertices $v_1$ and $v_2$ such that $|v_1|=1$ and $|v_2|=2$. Take $f=\sqrt{q+1}\chi_{v_1}$ and $g=\sqrt{q(q+1)}\chi_{v_2},$
where $\chi_v$ denotes characteristic function on the set $\{v\}$. Then it is easy to see that
$f, g \in \mathbb{T}_2$ with
$$\|f\|_2=\|g\|_2=\|f+g\|_2=\|f-g\|_2=1
$$
and hence the parallelogram law
$$\|f+g\|_2^2+\|f-g\|_2^2= 2(\|f\|_2^2+\|g\|_2^2)
$$
is not satisfied. Therefore, $\mathbb{T}_2$ cannot be a Hilbert space under $\|.\|_2$.
\brem
In the classical Hardy space $H^2$ of the unit disk,
$$\displaystyle\sup_{0\leq r<1}M_{2}(r,f)=\left(\frac{1}{2\pi}\int_{0}^{2\pi}|f(e^{i\theta})|^{2}\,d\theta\right)^{1/2},
$$
which is due to Littlewood's subordination theorem and mean convergence theorem (see \cite{Duren:Hpspace}).
Therefore $H^2$ becomes a Hilbert space in a natural way. On the other hand a similar situation does not occur in the
$\mathbb{T}_{p}$ spaces.
\erem

\bthm
For $f\in \mathbb{T}_{\infty}$, we have
$\lim\limits_{s\rightarrow\infty}\|f\|_s=\|f\|_\infty$.
\ethm
\bpf
For $n\in \mathbb{N}_{0}$ and $0<s<t\leq\infty$, we see that $M_s(n,f)\leq M_t(n,f)$  and thus,
$\|f\|_s\leq \|f\|_t$ for $s<t$ which in turn gives that
$\limsup\limits_{s\rightarrow\infty}\|f\|_s\leq\|f\|_\infty$.
On the other hand, for each $n\in \mathbb{N}_{0}$, we find that
$$
c_n^{-1/s} M_\infty(n,f)\leq M_s(n,f)\leq \|f\|_s,
$$
where $c_n$ is defined by \eqref{cn}.
Now, by letting $s\rightarrow\infty$ and taking supremum over $n\in \mathbb{N}_{0}$, we get
$\|f\|_\infty \leq \liminf\limits_{s\rightarrow\infty}\|f\|_s$.  Hence,
$\lim\limits_{s\rightarrow\infty}\|f\|_s=\|f\|_\infty$ as desired.
 \epf

\begin{Lem}\label{lem:bound} \emph{(Growth Estimate)}
{\rm (\cite[Lemma 3.12]{MP-Tp-spaces})}
Let $T$ be a $(q+1)$-homogeneous tree rooted at $\textsl{o}$ and $0<p<\infty$.
If $f$ is an element of $\mathbb{T}_{p}$ or $\mathbb{T}_{p,0}$, then we have
$$ |f(v)|\leq \{(q+1)q^{|v|-1}\}^{\frac{1}{p}} \|f\|_p ~\mbox{ for }~ v\in T.
$$
\end{Lem}

\begin{Lem}
{\rm (\cite[Theorems 3.10 and 3.11]{MP-Tp-spaces})}
For $0<p\leq\infty$, the space $\mathbb{T}_{p}$ is not  separable, whereas $\mathbb{T}_{p,0}$ is a
separable space as the span of $\{\chi_v :\, v\in T\}$ is
dense in $\mathbb{T}_{p,0}$.
\end{Lem}

%Proof of all these results stated in this section are straight forward from the definitions. For the proof, see \cite{MP-Tp-spaces}.

\section{Bounded Composition Operators} %\label{MP5Sec3}
\subsection{Bounded composition operators on $\mathbb{T}_{p}$} \label{MP5Subsec1}
A linear operator $A$ on a Banach space is said to be \textit{bounded}
if the operator norm $ \|A\|= \mbox{sup}\{\|Ax\| : \|x\|=1\}$ is finite. In this section, we  discuss
boundedness of composition operator $C_\phi$ on $\mathbb{T}_{p}$ spaces and compute their norm.

For the boundedness of $C_\phi$, we will discuss it case by case.

\begin{Thm}%\label{bdd:T_infty}
{\rm (\cite[Theorem 1]{CO-Tp-spaces})}
Every self-map $\phi$ of $T$ induces a bounded composition operator on $\T_\infty$ with $\|C_\phi\|=1$.
\end{Thm}

Next, we consider composition operators on $\T_p$ for $1\leq p<\infty$ over $2$-homogeneous trees. Every self-map $\phi$ of
$2$-homogeneous tree induces a bounded $C_\phi$ on $\mathbb{T}_{p}$ (see \cite{CO-Tp-spaces}).
\bthm\label{bdd:Tp,q=1}
Let $T$ be a $2$-homogeneous tree with root $\textsl{o}$ and let $D_n=\{a_n,b_n\}$ for each $n\in\mathbb{N}$.
Furthermore, let $\phi$ be a self-map of  $T$ and $C_\phi$ be the induced composition operator on $\mathbb{T}_{p}$
for $1\leq p<\infty$.
Then we have the following:
\begin{enumerate}
  \item If $\phi(\textsl{o})\neq \textsl{o}$, then $\|C_\phi\|^p=2$.
  \item If $\phi(\textsl{o})= \textsl{o}$, then  any one of the following distinct cases must occur:
  \begin{enumerate}
    \item Either $\phi \equiv \textsl{o}$ or for every $n\in\mathbb{N}$, if $\phi$ maps $D_n$ bijectively onto $D_m$ for some
    $m\in\mathbb{N}$  then $\|C_\phi\|^p=1$.
    \item If $\phi$ maps exactly one element of $D_n$ to $\textsl{o}$ for each $n\in\mathbb{N}$ then $\|C_\phi\|^p=\frac{3}{2}$.
    \item Either there exists an $n\in\mathbb{N}$ such that $\phi(a_n)=\phi(b_n)\neq \textsl{o}$ or if there exists an
    $n\in\mathbb{N}$ such that $|\phi(a_n)|$ and $|\phi(b_n)|$ are not equal and both are different from $0$ then $\|C_\phi\|^p=2$.
  \end{enumerate}
\end{enumerate}
\ethm
\bpf
From the growth estimate for $2$-homogeneous trees, it follows that for each $n\in\mathbb{N}_0$,
$$ M_p^p(n,C_\phi f) = \frac{1}{c_n}\sum\limits_{|v|=n} |f(\phi(v))|^p \leq 2\|f\|^p ~\mbox{ for every }~ f\in \mathbb{T}_p.
$$
This yields that $\|C_\phi\|^p\leq 2$. Thus every self-map $\phi$ of $T$ induces a
bounded $C_\phi$ on $\mathbb{T}_p$ with $\|C_\phi\|^p\leq 2$.

Suppose that $w=\phi(\textsl{o})\neq \textsl{o}$. For $f=\displaystyle{ 2^{\frac{1}{p}}\chi_{w}}$, we have $\|f\|=1$ and
$\|C_\phi(f)\|^p=2$ and hence $\|C_\phi\|^p=2$.

Now suppose that $\phi(\textsl{o})= \textsl{o}$. Then we need to consider all the five possible cases.

Suppose that $\phi \equiv \textsl{o}$. Then for each $n\in\mathbb{N}_0$,
$$M_p^p(n,C_\phi f)= |f(\textsl{o})|^p \leq \|f\|^p ~\mbox{ for every }~ f\in \mathbb{T}_p.
$$
This yields that $\|C_\phi\|^p\leq1$.
For $f=\chi_{\textsl{o}}$, we obtain that $\|f\|^p=\|C_\phi(f)\|^p$ and thus, $\|C_\phi\|^p=1$.

Suppose that for every $n\in\mathbb{N}$, $\phi$ maps $D_n$ bijectively onto $D_m$ for some $m\in\mathbb{N}$.
Then,   $M_p^p(n,C_\phi f)=M_p^p(m,f)$ for every $n\in\mathbb{N}$ and for some $m\in\mathbb{N}$.
Thus
$$\|C_\phi f\|^p\leq \|f\|^p ~\mbox{ for every }~ f\in \mathbb{T}_p,
$$
which gives that $\|C_\phi\|^p\leq 1$.
As in the previous case, by considering $f=\chi_{\textsl{o}}$,   we get $\|C_\phi\|^p=1$.

Suppose that $\phi$ maps exactly one element of $D_n$ to $\textsl{o}$ for each $n\in\mathbb{N}$. Then,
in view of growth estimate for $2$-homogeneous trees along with this assumption, we see that
$$\|C_\phi f\|^p\leq\frac{3}{2}\|f\|^p ~\mbox{ for every }~ f\in \mathbb{T}_p
$$
which gives  $\|C_\phi\|^p\leq 3/2.$  On the other hand, by assumption, either $a_1$ or $b_1$ maps to $\textsl{o}$.
Without loss of generality, we assume that $\phi(a_1)=\textsl{o}$. Take $\phi(b_1)=w$ and
$f=\displaystyle{\chi_{\textsl{o}}+ 2^{\frac{1}{p}}\chi_{w}}$. Then, $\|f\|=1$ and
$$M_p^p(1,C_\phi f)=\frac{3}{2}=\|C_\phi(f)\|^p.
$$
Thus, $\|C_\phi\|^p=3/2$.

Now assume that there exists an $n\in\mathbb{N}$ such that $w=\phi(a_n)=\phi(b_n)\neq \textsl{o}$.
We have already observed that $\|C_\phi\|^p\leq 2$. For $f=\displaystyle{ 2^{\frac{1}{p}}\chi_{w}}$, we have
$$\|f\|=1 ~\mbox{ and }~ \|C_\phi(f)\|^p=2
$$
and therefore, $\|C_\phi\|^p=2$.

Finally, assume that there exists an $n\in\mathbb{N}$ such that $|\phi(a_n)|$ and $|\phi(b_n)|$ are not equal
and are different from $0$. Now, we take
$$f=\displaystyle{ 2^{\frac{1}{p}}(\chi_{u}+\chi_{v})},
$$
where $\phi(a_n)=u$ and $\phi(b_n)=v$. It follows that $\|f\|=1$ and $\|C_\phi(f)\|^p=2$, which gives that  $\|C_\phi\|^p=2$.
\epf

\bdefn
A bijective self-map $\phi$ of $T$ is called  an \textit{automorphism} of $T$ if any two vertices $v$ and $w$ are
neighbours if and only if $\phi(v)$ and $\phi(w)$ are neighbours.
The set of all automorphisms of $T$ is denoted by  $\Aut (T)$.
\edefn
From Theorem \ref{bdd:Tp,q=1} and \cite[Theorems 1 and 5]{CO-Tp-spaces}, we have the following result:
\bcor
Let $C_\phi$ be  a composition operator on $\mathbb{T}_{p}$ induced by an automorphic symbol $\phi$
of $T$. Then we have the following:
\begin{enumerate}
\item[(i)]  $\|C_\phi\| = 1$ if $p=\infty$.
\item[(ii)] For $q\geq 1$ and $1\leq p <\infty$,
$$ \|C_\phi\|^p= \left \{\begin{array}{ll}
(q+1) q^{|\phi(\textsl{o})|-1}   & ~\mbox{ if }~\phi(\textsl{o})\neq \textsl{o}, \\
1 &~\mbox{ if }~ \phi(\textsl{o})= \textsl{o}.
\end{array} \right.
$$
\end{enumerate}
\ecor

Next, we consider  composition operators on $\T_p$ for $1\leq p<\infty$ over $(q+1)$-homogeneous tree with $q\geq 2$.
A self-map $\phi$ of $T$ is called  \textit{bounded} if $\{|\phi(v)| :v\in T\}$ is a bounded set in $\mathbb{N}_0$.
From \cite[Theorem 3]{CO-Tp-spaces},
it is easy to see that every bounded self-map of $T$ induces bounded composition operators.
\bthm%\label{}
If $T$ is a $(q+1)$-homogeneous tree with $q\geq 2$ and $\phi$ is a bounded self-map of $T$
such that $\sup\limits_{v\in T} |\phi(v)| =M$,
then $\|C_\phi\|^p \leq c_M.$ Moreover, $\|C_\phi\|^p=c_M$ if and only if
$$\displaystyle{\sup\limits_{n\in \mathbb{N}_0} \frac{N _{M,n}}{c_n}=1 }.
$$
\ethm
\bpf
For $n\in \mathbb{N}_0$ and $f\in \mathbb{T}_p$,   by  Lemma \Ref{lem:bound}, we have
$$M_p^p(n,C_\phi f)\leq c_M \|f\|^p.
$$
Thus, $\phi$ induces a bounded $C_\phi$ with $\|C_\phi\|^p \leq c_M.$

Let us now prove the equality case. Suppose that
$\displaystyle{\sup\limits_{n\in \mathbb{N}_0} \frac{N _{M,n}}{c_n}=1 }$. Then there are two cases. First
we consider the case  $N_{M,k}=c_k$ for some $k\in \mathbb{N}_0$. This
means that $\phi: D_k\rightarrow D_M$ is a constant, say, $\phi(v)=w \in D_M$ for all $v\in D_k$. For $f= (c_M)^{\frac{1}{p}}\chi_{w}$, we have
$$\|f\|=1 ~\mbox{ and }~\|C_\phi(f)\|^p=c_M %\|f\|^p
$$
which proves that $\|C_\phi\|^p=c_M$.

Next, we suppose that $N _{M,k}\neq c_k$ for all $k\in \mathbb{N}_0$. Then, there is a sequence $\{n_k\}$ such that
$$\displaystyle \frac{N _{M,n_k}}{c_{n_k}}\rightarrow 1 ~\mbox{ as $k\rightarrow \infty$}.
$$
For each $k\in \mathbb{N}$, choose $w_k \in D_M$ such that $N_{M,n_k}=N_{\phi}(n_k,w_k)$. Take $f_k= (c_M)^{\frac{1}{p}}\chi_{w_k}$ so
that $\|f_k\|=1$ and
$$c_M \frac{N_{M,n_k}}{c_{n_k}} \leq M_p^p(n_k, C_\phi(f)) \leq \|C_\phi(f)\|^p \leq \|C_\phi\|^p.
$$
By allowing $k\rightarrow\infty$, we get $c_M \leq \|C_\phi\|^p$, and thus, $\|C_\phi\|^p = c_M$ in either case.

For the converse part, we assume that $\|C_\phi\|^p=c_M$.  Suppose  on the contrary that $\displaystyle{\sup\limits_{n\in \mathbb{N}_0} \frac{N _{M,n}}{c_n} \leq \delta <1}$.
Then, $N _{M,n} \leq \delta c_n$ for every $n$. Note that, there are at least $N _{M,n}$ vertices from $D_n$
mapped into $D_M$ and therefore, there are at most $c_n-N _{M,n}$ vertices of $D_n$ mapped into $\{v:\, |v|<M\}$ for each $n$.
For $f\in \mathbb{T}_{p}$, we obtain
\beqq
M_p^p(n,C_\phi f)&=&\frac{1}{c_n}\sum_{\substack{|\phi(v)|=M \\ |v|=n}} |f(\phi (v))|^p +
\frac{1}{c_n}\sum_{\substack{|\phi(v)|<M \\ |v|=n}} |f(\phi (v))|^p \\
&\leq & \frac{N_{M,n}}{c_n} c_M \|f\|^p+ \frac{c_n-N_{M,n}}{c_n} c_{M-1} \|f\|^p \\
&\leq & \{ 1+(q-1)\delta \} c_{M-1} \|f\|^p.
\eeqq
Therefore, $\|C_\phi\|^p \leq \{ 1+(q-1)\delta \} c_{M-1} < c_M$, which is a contradiction. Hence,
$\|C_\phi\|^p=c_M$ if and only if
$\displaystyle{\sup\limits_{n\in \mathbb{N}_0} \frac{N _{M,n}}{c_n}=1 }$.
\epf

Now, we consider   general self-maps on $(q+1)$-homogeneous trees.
\bthm\label{norm}
Let $T$ be a  $(q+1)$-homogeneous tree and  $1\leq p<\infty$. Then $C_\phi$ is bounded on $\mathbb{T}_{p}$ if and only if
$$ \alpha:= \ds \sup\limits_{n\in \mathbb{N}_0} \left \{ \frac{1}{c_n}\sum\limits_{m=0}^{\infty}N_{m,n}c_m \right\}< \infty.
$$
Moreover, $\|C_\phi\|^p=\alpha$.
\ethm
\bpf
Assume that $\alpha < \infty.$ First we show that $C_\phi$ is bounded on $\mathbb{T}_{p}$. To do this,
for $n\in \mathbb{N}_0$ and $f\in \mathbb{T}_{p}$, we find that
\beqq
M_p^p(n,C_\phi f)
&=& \frac{1}{c_n} \left \{ \sum\limits_{m=0}^{\infty} \sum_{\substack{|\phi(v)|=m \\ |v|=n}} |f(\phi (v))|^p \right \}\\
&\leq& \left \{ \frac{1}{c_n}\sum\limits_{m=0}^{\infty}N_{m,n}c_m \right\} \|f\|^p\\
&\leq& \alpha  \|f\|^p,
\eeqq
which yields that $C_\phi$ is bounded on $\mathbb{T}_{p}$ and
\be\label{mp5-eq3}
\|C_\phi\|^p\leq \alpha.
\ee
%=\ds \sup\limits_{n\in \mathbb{N}_0} \left \{ \frac{1}{c_n}\sum\limits_{m=0}^{\infty}N_{m,n}c_m \right\}
Conversely, suppose that $C_\phi$ is bounded on $\mathbb{T}_{p}$. In order to show that $\alpha$ is finite,
we fix $n \in \mathbb{N}_0$. For each $m \in \mathbb{N}_0$, choose $v_m\in D_m$ such that
$N_\phi(n,v_m)=N_{m,n}$. Take $f=\sum\limits_{m=0}^\infty (c_m)^{\frac{1}{p}}\chi_{v_m}$,
so that $\|f\|=1$ and
$$M_p^p(n,C_\phi f)= \ds \frac{1}{c_n}\sum\limits_{m=0}^{\infty}N_{m,n}c_m ,$$
which gives that
\be \label{mp5-eq4}
\alpha=\ds \sup\limits_{n\in \mathbb{N}_0} \left \{ \frac{1}{c_n}\sum\limits_{m=0}^{\infty}N_{m,n}c_m \right\} \leq \|C_\phi\|^p,
\ee
and hence the desired result follows. Moreover, by \eqref{mp5-eq3} and \eqref{mp5-eq4}, it follows that $\|C_\phi\|^p=\alpha$.
\epf

\subsection{Special symbols}\label{MP5Subsec2}
Every self-map $\phi$ of $T$ induces a bounded operator $C_\phi$ on $\mathbb{T}_\infty$, or $\mathbb{T}_p$ spaces over $2$-homogeneous trees.
Unlike the classical Hardy space settings, there are self-maps $\phi$ of $T$ which do not induce bounded $C_\phi$ on
$\mathbb{T}_p$ with $1\leq p<\infty$ over $(q+1)$-homogeneous trees, $q\geq 2$ (See \cite[Section 5]{CO-Tp-spaces}).

The following example shows that there are bijective self-maps of $T$ which do not induce bounded composition
operator $C_\phi$ for $(q+1)$-homogeneous trees with $q\geq 2$.
\beg\label{mp5-eg-1}
For each $n\in \N$ which is not of the form $n=4k$, $k\in \N_0$, choose $v_n\in T$ such that $|v_n|=n$. Define
$$ \phi(v)= \left \{\begin{array}{ll}
v_{4k+2} & ~\mbox{ if }~ v=v_{2k+1}  ~\mbox{ for some }~ k\in \N_0,\\
v_{2k+1} &~\mbox{ if }~ v=v_{4k+2}~\mbox{ for some }~ k\in \N,\\
v & ~\mbox{ elsewhere.}
\end{array} \right.
$$
Clearly, $\phi$ is bijective on $T$. For $k\in \N$, let
$f_k=\ds (c_{4k+2})^{\frac{1}{p}}\chi_{v_{4k+2}}$. Then $\|f\|=1$ and
$$ \|C_\phi\|^p \geq \|C_\phi(f_k)\|^p \geq M_p^p(2k+1,C_\phi(f_k))=q^{2k+1}.
$$
Since $q\geq 2$, it follows that $C_\phi$ is an unbounded operator on $\mathbb{T}_p$.
\eeg
Motivated by the above example, we wish to characterize all the bounded composition operators that are induced by
univalent (injective) symbols (see Corollary \ref{mp5-inj.sym2}).

\bthm\label{mp5-inj.sym1}
Let $\phi$ be a self-map of $(q+1)$-homogeneous tree $T$ with $q\geq2$, and $1\leq p<\infty$.
If $C_\phi$ is bounded on $\mathbb{T}_p$, then there exists an $M>0$ such that $|\phi(v)|\leq |v|+M$ for
all $v\in T$.
\ethm
\bpf
Suppose that $C_\phi$ is bounded on $\mathbb{T}_p$. Set $a_n=\max\limits_{|v|=n} |\phi(v)|$ for $n\in \N_0$, and
for each $n$, choose $v_n\in D_n$ such that
$|\phi(v_n)|=a_n$. Furthermore, for each $n$, take
$f_n=\ds (c_{a_n})^{\frac{1}{p}}\chi_{\phi(v_n)}$. Then
$$
M_p^p(n,C_\phi f_n)=q^{a_n-n}\leq \|C_\phi\|^p,
$$
which gives that $\{a_n-n\}$ is a bounded sequence. The desired result follows.
\epf

Converse of Theorem \ref{mp5-inj.sym1} holds if, in addition, $\phi$ is injective or finite-valent.
\bcor\label{mp5-inj.sym2}
If $\phi$ is an injective self-map of $(q+1)$-homogeneous tree $T$ with $q\geq2$ and $1\leq p<\infty$, then
$C_\phi$ is bounded on $\mathbb{T}_p$ if and only if  there exists an $M>0$ such that $|\phi(v)|\leq |v|+M$ for
all $v\in T$.
\ecor
\bpf
Suppose that there exists an $M>0$ such that $|\phi(v)|\leq |v|+M$ for
all $v\in T$. Therefore, $a_n\leq n+M$ for all $n$, where $a_n$ is taken as in Theorem \ref{mp5-inj.sym1}.
For an arbitrary function $f$ with $\|f\|=1$, we have
\beqq
M_p^p(n,C_\phi f)
&\leq& \frac{1}{c_n}  \sum\limits_{m=0}^{a_n} \sum\limits_{|w|=m} |f(w)|^p ~\mbox{ (since $\phi$ is injective)}\\
&\leq& \frac{1}{c_n} \sum\limits_{m=0}^{a_n}c_m
= \frac{1}{c_n}+ \frac{1}{(q+1)\,q^{n-1}}\sum\limits_{m=1}^{a_n}(q+1)\,q^{m-1}\\
%&=&\frac{1}{c_n}+ \frac{1}{q^{n-1}}\sum\limits_{k=0}^{a_n-1}q^{k}\\
&=&\frac{1}{c_n}+\frac{q^{a_n}-1}{q^{n-1}(q-1)}\\
&\leq& %\frac{1}{c_n}+\frac{q^{a_n}}{q^{n-1}(q-1)}=
\frac{1}{c_n}+ q^{a_n-n}\frac{q}{q-1}\\
%=\frac{1}{c_n}+ q^{a_n-n}\frac{q}{q-1}
&\leq& 1+ q^M\left(\frac{q}{q-1}\right).
\eeqq

Thus, $C_\phi$ is bounded on $\mathbb{T}_p$. The converse part is a consequence of Theorem~\ref{mp5-inj.sym1}.
\epf

\bdefn
Let $\phi$ be a self-map of $T$ and $k\in \N$ be fixed. We say that $\phi$ is $k$-valent map if every vertex of $T$
has at most $k$ pre-images and there is a vertex of $T$ which has exactly $k$ pre-images. The map $\phi$ is said to be
\textit{finite-valent} if there exists an $k\in \N$ such that $\phi$ is $k$-valent.
\edefn
The next corollary follows directly from the proofs of Theorem \ref{mp5-inj.sym1} and Corollary \ref{mp5-inj.sym2}.
\bcor\label{mp5-finite-valant}
%Let $T$ and $p$ be as in Theorem \ref{mp5-inj.sym} and $\phi$ be a finite-valent map. Then
Let $\phi$ be an finite-valent self-map of $(q+1)$-homogeneous tree $T$ with $q\geq2$, and $1\leq p<\infty$.
$C_\phi$ is bounded on $\mathbb{T}_p$ if and only if  there exists an $M>0$ such that $|\phi(v)|\leq |v|+M$ for
all $v\in T$.
\ecor
\brem
The assumption that $\phi$ is finite-valent is necessary in Corollary \ref{mp5-finite-valant}. To see this, for each $n$, fix $v_n\in D_n$ and $\phi(v)=v_n$ if $|v|=n$.
For each $n$, choose $f_n=(c_n)^{\frac{1}{p}}\chi_{v_n}$ so that
$$M_p^p(n,C_\phi f_n)=(q+1)q^{n-1}\leq \|C_\phi\|^p,
$$
which gives that $C_\phi$ cannot be a bounded operator.
\erem

\subsection{Composition operators on $\mathbb{T}_{p,0}$}\label{MP5Subsec3}
\bthm\label{little Tp}
Let $\phi$ be a self-map of $T$ and $f\in \mathbb{T}_{p}$. If $|\phi(v)|\rightarrow\infty$
 and $|f(v)|\rightarrow0$ as
$|v|\rightarrow\infty$, then $|C_\phi f(v)|\rightarrow0$ as $|v|\rightarrow\infty$.
\ethm
\bpf
Assume the hypothesis. Let $\epsilon>0$ be given. Then there exists an $N_1\in\N$ such that
$|f(w)|<\epsilon$  for all $|w|\geq N_1$. Given $N_1>0$, there exists an $N\in\N$ such that
$|\phi(v)|\geq N_1$  for all $|v|\geq N$. This gives that
$|C_\phi f(v)|<\epsilon$  for all $|v|\geq N$, i.e., $|C_\phi f(v)|\rightarrow0$ as $|v|\rightarrow\infty$.
\epf

\blem\label{rem:Tp0}
\begin{enumerate}
\item $h\in\mathbb{T}_{\infty,0}$ if and only if $|h(v)|\rightarrow0$ as $|v|\rightarrow\infty$.
\item Let $T$ be a $2$-homogeneous tree. For $1\leq p<\infty$,
 $h\in\mathbb{T}_{p,0}$ if and only if $|h(v)|\rightarrow0$ as $|v|\rightarrow\infty$.
\end{enumerate}
\elem

\bpf
 For $n\in \N$ and $h\in\mathbb{T}_\infty$, we have
$$M_\infty(n,h)=\max\limits_{|v|=n}|h(v)|=h(v_n) ~\mbox{ for some $v_n$ with $|v_n|=n$}.
$$
In view of this, it is easy to see that $h\in\mathbb{T}_{\infty,0}$ if and only if $M_\infty(n,h)\rightarrow0$
as $n\rightarrow\infty$ if and only if $|h(v)|\rightarrow0$ as $|v|\rightarrow\infty$.

Let $T$ be a $2$-homogeneous tree and for $n\in \N$, take $D_n=\{a_n,b_n\}$. For $h\in\mathbb{T}_p$, we have
$$M_p^p(n,h)=\frac{1}{2}(|h(a_n)|^p+|h(b_n)|^p).
$$
This yields that $h\in\mathbb{T}_{p,0}$ if and only if $M_p(n,h)\rightarrow0$
as $n\rightarrow\infty$ if and only if $|h(v)|\rightarrow0$ as $|v|\rightarrow\infty$.
\epf

\bthm\label{bdd:T_infty,0}
$C_\phi$ is bounded operator on $\mathbb{T}_{\infty,0}$ if and only if
$|\phi(v)|\rightarrow\infty$ as $|v|\rightarrow\infty$.
\ethm
\bpf
Since $\|C_\phi(f)\|_\infty\leq\|f\|_\infty$ for all $f\in T_\infty$, it is enough to prove that
$\mathbb{T}_{\infty,0}$ is invariant under $C_\phi$ if and only if $|\phi(v)|\rightarrow\infty$ as
$|v|\rightarrow\infty$.

Suppose that $|\phi(v)|\rightarrow\infty$ as $|v|\rightarrow\infty$. Then by  Theorem
\ref{little Tp} and Lemma \ref{rem:Tp0}, $|C_\phi f(v)|\rightarrow0$ as $|v|\rightarrow\infty$ for all $f\in T_{\infty,0}$.
That is, $\mathbb{T}_{\infty,0}$ is invariant under $C_\phi$.

For the converse part, assume that $|\phi(v)|\not\rightarrow\infty$ as $|v|\rightarrow\infty$.
Then there exists a sequence $\{v_k\}$ and $M>0$ such that $|v_k|\geq k$ and $|\phi(v_k)|\leq M$
for all $k\in \N$. Define
$$ f(v)= \left \{\begin{array}{ll}
1 & \mbox{ if }~ |v|\leq M,\\
1/|v_k| & \mbox{ if $|v|>M$ and $v=v_k$ for some }~ k\in \N,\\
0 & \mbox{  elsewhere.}
\end{array} \right.
$$
Then $f\in \mathbb{T}_{\infty,0}$. But $M_\infty(f\circ \phi, |v_k|)=1$ for all $k$ and so,
$f\circ \phi \notin \mathbb{T}_{\infty,0}$. This completes the proof.
\epf

\bthm\label{bdd:Tp,0,q=1}
Let $T$ be a $2$-homogeneous tree with root $\textsl{o}$ and let $D_n=\{a_n,b_n\}$ for each $n\in\mathbb{N}$ and $\phi$ be a
self-map of  $T$. Then, $C_\phi$ is a bounded operator on $\mathbb{T}_{p,0}$, $1\leq p<\infty$, if and only if
$|\phi(v)|\rightarrow\infty$ as $|v|\rightarrow\infty$. Moreover, we have the following norm estimates.
\begin{enumerate}
  \item If $\phi(\textsl{o})\neq \textsl{o}$, then $\|C_\phi\|^p=2$.
  \item If $\phi(\textsl{o})= \textsl{o}$, then  any one of the following distinct cases must occur:
  \begin{enumerate}
   \item For every $n\in\mathbb{N}$, if $\phi$ maps $D_n$ bijectively onto $D_m$ for some
    $m\in\mathbb{N}$,  then $\|C_\phi\|^p=1$.
  \item Either there exists an $n\in\mathbb{N}$ such that $\phi(a_n)=\phi(b_n)\neq \textsl{o}$ or if there exists an
    $n\in\mathbb{N}$ such that $|\phi(a_n)|$ and $|\phi(b_n)|$ are not equal and different from $0$ then $\|C_\phi\|^p=2$.
  \end{enumerate}
\end{enumerate}
\ethm
\bpf
For $2$-homogeneous  trees, every self-map $\phi$ induces bounded composition operators on $\mathbb{T}_{p}$.
Therefore it suffices to prove that  $\mathbb{T}_{p,0}$ is invariant under $C_\phi$ if and only if
$|\phi(v)|\rightarrow\infty$ as $|v|\rightarrow\infty$.

Necessary part follows from Theorem \ref{little Tp} and Lemma \ref{rem:Tp0}. For the proof of the sufficiency part,
assume on the contrary that $|\phi(v)|\not\rightarrow\infty$ as $|v|\rightarrow\infty$.
Take $f$ as in Theorem \ref{bdd:T_infty,0}. Then $f\in \mathbb{T}_{p,0}$. But $f(\phi(v_k))=1$ for all
$k$ which gives $M_p^p(f\circ \phi, |v_k|)\geq 1/2$  for all $k$ and so,
$f\circ \phi \notin \mathbb{T}_{p,0}$. This completes the proof.

The proof of norm estimates is similar to that of Theorem \ref{bdd:Tp,q=1}.
\epf

\bthm%\label{norm}
Let  $q\geq 2$ and  $1\leq p<\infty$. If $\frac{1}{c_n}\sum\limits_{m=0}^{\infty}N_{m,n}c_m\rightarrow0$
as $n\rightarrow\infty$, then $C_\phi$ is bounded on $\mathbb{T}_{p,0}$. Moreover, $\|C_\phi\|^p=\alpha$,
where
$$ \alpha=\ds \sup\limits_{n\in \mathbb{N}_0} \left \{ \frac{1}{c_n}\sum\limits_{m=0}^{\infty}N_{m,n}c_m \right\}.
$$
\ethm
\bpf
By Theorem \ref{norm}, for the boundedness of $C_\phi$ on $\mathbb{T}_{p,0}$, it is enough to prove that
$C_\phi$ maps $\mathbb{T}_{p,0}$ into $\mathbb{T}_{p,0}$. It follows from the proof of Theorem \ref{norm} that
$$M_p^p(n,C_\phi f) \leq \left \{ \frac{1}{c_n}\sum\limits_{m=0}^{\infty}N_{m,n}c_m \right\} \|f\|^p,
$$
which forces that  $\mathbb{T}_{p,0}$ is invariant under $C_\phi$  whenever $\frac{1}{c_n}\sum\limits_{m=0}^{\infty}N_{m,n}c_m\rightarrow0$
as $n\rightarrow\infty$. Moreover, we have $\|C_\phi\|^p\leq \alpha$.

To prove that the equality holds in the last inequality, we fix $n \in \mathbb{N}_0$.
For each $m \in \mathbb{N}_0$, choose $v_m\in D_m$ such that $N_\phi(n,v_m)=N_{m,n}$. Then,
$$f_k=\sum\limits_{m=0}^k(c_m)^{\frac{1}{p}}\chi_{v_m} \in \mathbb{T}_{p,0} ~\mbox{ and $\|f_k\|=1$ for all $k$.}
$$
Therefore,
$$  \frac{1}{c_n}\sum\limits_{m=0}^k N_{m,n}c_m =M_p^p(n,C_\phi f_k)\leq \|C_\phi\|^p ~\mbox{ for all}~ k,
$$ which gives that $\alpha \leq \|C_\phi\|^p$,
and this completes the proof.
\epf

\bthm
If $C_\phi$ is bounded on $\mathbb{T}_{p,0}$, $1\leq p <\infty$, then
\be\label{MP5-eq1}
\frac{c_{|v|}}{c_n} N_{\phi}(n,v)\ra 0 ~\mbox{ as }~ n\ra \infty ~\mbox{ for every }~ v\in T.
\ee
%where $N_{\phi}(n,v)$ is defined in Subsection \ref{MP5Subsec1}.
\ethm
\bpf
For each $v\in T$, define $f_v=  (c_{|v|})^{\frac{1}{p}} \chi_v$.
Then $f_v \in \mathbb{T}_{p,0}$ with $\|f_v\|=1$ and
 %for each fixed $v\in T$, we have for $n\in \mathbb{N}_0$,
$$
M_p^p(n,C_\phi f_v)= \frac{c_{|v|}}{c_n} N_{\phi}(n,v).
$$
Since $C_\phi(\mathbb{T}_{p,0})\subseteq \mathbb{T}_{p,0}$, we have
$\frac{c_{|v|}}{c_n} N_{\phi}(n,v)\ra 0 ~\mbox{ as }~ n\ra \infty$ for every $v\in T$.
\epf
\brem
The condition \eqref{MP5-eq1} is equivalent to $C_\phi(\chi_v)\in \mathbb{T}_{p,0}$ for every
$v\in T$. This in turn is also equivalent to saying that $C_\phi(E)\subseteq \mathbb{T}_{p,0}$, where
$E=~\mbox{ Span}\{ \chi_v: v\in T\}$ is a dense subspace of $\mathbb{T}_{p,0}$ under $\|.\|_p$.
\erem

\section{Invertible Composition Operators}\label{MP5Sec4}

Recall that a bounded linear operator $A$ on a normed linear space $X$ is said to be an \textit{invertible}
if there exists a bounded linear operator $B$ on $X$ such that $B(A(x))=A(B(x))=x$ for all $x\in X$.
Such an operator $B$ is called inverse of $A$ and is denoted by $A^{-1}$.
\blem\label{mp5-inv-lem}
If $C_\phi$ is an invertible operator on $\mathbb{T}_p$, $p\geq 1$, then $\phi$ is bijective on $T$.
Moreover, $C_\phi^{-1}=C_{\phi^{-1}}$.
\elem
\bpf
Assume that $C_\phi$ is invertible.

Suppose on the contrary that $\phi$ is not onto. Pick a vertex $w\in T\setminus \phi(T)$, where $\phi(T)$ denotes the
image of $T$ under $\phi$. Then for $f=\chi_w$, $f\not\equiv0$
and $C_\phi(f)=0$. Therefore, $C_\phi$ is not injective which leads to a contradiction. Hence $\phi$ is onto.

Suppose on the contrary that $\phi$ is not injective on $T$. Then there exists $v_1,v_2\in T$ such that $v_1\neq v_2$ and
$\phi(v_1)=\phi(v_2)=w$. Take $g=\chi_{v_1}\in \mathbb{T}_p$. But there is no $f\in \mathbb{T}_p$
such that $C_\phi(f)=g$, because $0=g(v_2)=f(w)=g(v_1)=1$. Therefore, $C_\phi$ is not onto which is again a
contradiction. Thus $\phi$ is injective.

Since $C_\phi$ is invertible, $\phi$ is bijective and there is a bounded linear operator $S$ on $\mathbb{T}_p$ such that
$C_\phi \circ S=S\circ C_\phi=I$, where $I$ is the identity operator on $\mathbb{T}_p$. Now, it is easy to see that
$S=C_{\phi^{-1}}$. The desired conclusion follows.
\epf

\bthm\label{mp5-inv-rem}
A bounded operator $C_\phi$ on $\mathbb{T}_p$ is invertible if and only if $\phi$ is bijective on $T$ and
$C_{\phi^{-1}}$ is a bounded operator on $\mathbb{T}_p$.
\ethm
\bpf
Suppose $C_\phi$ is an invertible operator on $\mathbb{T}_p$. Then by Lemma \ref{mp5-inv-lem}, $\phi$ is bijective
on $T$ and $C_\phi^{-1}=C_{\phi^{-1}}$ is a bounded operator on $\mathbb{T}_p$. Converse holds trivially, since
$C_{\phi^{-1}}$ will be an inverse of $C_\phi$.
\epf

Since every self-map $\phi$ of $T$ induces a bounded operator $C_\phi$ on $\mathbb{T}_\infty$  (resp.  on $\mathbb{T}_p$
spaces over $2$-homogeneous trees), it is easy to obtain the following results.
\bcor
The operator $C_\phi$ is invertible on $\mathbb{T}_\infty$ if and only if $\phi$ is bijective on $T$.
\ecor
\bcor
Let $T$ be a $2$-homogeneous tree and let $1\leq p< \infty$.
The operator $C_\phi$ is invertible on $\mathbb{T}_p$ if and only if $\phi$ is bijective on $T$.
\ecor
\bcor
\begin{enumerate}
\item The operator $C_\phi$ is invertible on $\mathbb{T}_{\infty,0}$ if and only if $\phi$ is bijective on $T$
and $|\phi(v)|\rightarrow\infty$ and $|\phi^{-1}(v)|\rightarrow\infty$ as $|v|\rightarrow\infty$.
\item The operator $C_\phi$ is invertible on $\mathbb{T}_{p,0}$ space over $2$-homogeneous trees if and only if
$\phi$ is bijective on $T$ and $|\phi(v)|\rightarrow\infty$ and $|\phi^{-1}(v)|\rightarrow\infty$ as $|v|\rightarrow\infty$.
\end{enumerate}
\ecor
Example \ref{mp5-eg-1} shows that there are bijective self-maps of $T$ which do not induce
bounded composition operator $C_\phi$ for the case of $(q+1)$-homogeneous trees with $q\geq 2$. Indeed,
there are bijective self-maps $\phi$ of $T$ which induce a bounded composition operator $C_\phi$
on $\mathbb{T}_p$ over $(q+1)$-homogeneous trees with $q\geq 2$, but $\phi^{-1}$ does not necessarily induce
a bounded composition operator $C_{\phi^{-1}}$.
\beg
For each $k\in \N$, choose a subset $A_{2k-1}$ of $k-1$ vertices in $D_{2k-1}$ and choose a
subset $A_{2k}$ of $k$ vertices in $D_{2k}$. Label the elements of $A_n$ as
$A_n=\{v_{n,1},v_{n,2},v_{n,3},\ldots\}$ for  each $n\in \N$.

Define $\phi$ as follows: $\phi(\textsl{o})=\textsl{o}$, $\phi(v)=v$ if $v\in D_k\setminus A_k$
and $\phi(v_{2k,1})=v_{k,1}$. For each $k\in \N$, $A_{2k-1}$ and $A_{2k}\setminus\{v_{2k,1}\}$ have the
same number of elements. We can thus define $\phi: A_{2k-1} \rightarrow A_{2k}\setminus\{v_{2k,1}\}$ bijectively
and so does for defining $\phi: A_{2k}\setminus\{v_{2k,1}\} \rightarrow A_{2k+1}\setminus\{v_{2k+1,1}\}$
bijectively. Thus, $\phi:T\rightarrow T$ becomes a bijective self-map of $T$.

Take an arbitrary function $f\in \mathbb{T}_p$ with $\|f\|=1$. Fix $n=2k-1$ for some $k\in \N$. Then
\beqq
M_p^p(n,C_\phi f)
&=& \frac{1}{c_n} \left( \sum\limits_{v\in D_n\setminus A_n} |f(\phi(v))|^p +
\sum\limits_{v\in A_n} |f(\phi(v))|^p \right)\\
&=& \frac{1}{c_n} \left( \sum\limits_{w\in D_n\setminus A_n} |f(w)|^p +
\sum\limits_{w\in A_{n+1}\setminus\{v_{n+1,1}\}} |f(w)|^p \right)\\
&\leq& \frac{1}{c_n} \left(c_n\|f\|^p+c_{n+1} \|f\|^p \right) =(1+q)\|f\|^p=1+q.
\eeqq
Next, fix $n=2k$ for some $k\in \N$. Then
\beqq
M_p^p(n,C_\phi f)
&=& \frac{1}{c_n} \left( \sum\limits_{v\in D_n\setminus A_n} |f(\phi(v))|^p +
\sum\limits_{v\in A_n\setminus\{v_{n,1}\}} |f(\phi(v))|^p+ |f(\phi(v_{n,1}))|^p \right)\\
&=& \frac{1}{c_n} \left( \sum\limits_{w\in D_n\setminus A_n} |f(w)|^p +
\sum\limits_{w\in A_{n+1}\setminus\{v_{n+1,1}\}} |f(w)|^p + |f(v_{k,1})|^p\right)\\
&\leq& \frac{1}{c_n} \left(c_n\|f\|^p+c_{n+1} \|f\|^p+c_k\|f\|^p \right) \leq(2+q)\|f\|^p=2+q.
\eeqq
Thus, $\phi$ induces a bounded composition operator with $\|C_\phi\|^p\leq 2+q$.

Finally, we consider the composition operator induced by $\phi^{-1}$. Recall that
$\phi^{-1}(v_{n,1})=v_{2n,1}$ for each $n$. For $n\in \N$, take
$f_n=\ds (c_{2n})^{\frac{1}{p}}\ds\chi_{v_{2n,1}}$. Then
$$
M_p^p(n,C_{\phi^{-1}} f_n)=q^{n}\leq \|C_{\phi^{-1}}\|^p,
$$
which gives that $\phi^{-1}$ cannot induce a bounded composition operator.
\eeg

The following result characterizes the invertible composition operators on $\mathbb{T}_p$ over
$(q+1)$-homogeneous trees  with $q\geq 2$.
\bthm
Let $T$ be a $(q+1)$-homogeneous tree with $q\geq 2$, and $1\leq p< \infty$. $C_\phi$ is invertible on
$\mathbb{T}_p$,  if and only if $\phi$
is bijective and there exists an $M>0$ such that $|\,|\phi(v)|-|v|\,|\leq M$ for all $v\in T$.
\ethm
\bpf
Suppose $C_\phi$ is an invertible operator on
$\mathbb{T}_p$. Then, by Theorem \ref{mp5-inv-rem}, $\phi$ is bijective on $T$ and both $C_\phi$, $C_{\phi^{-1}}$ are bounded
 operators on $\mathbb{T}_p$. By Corollary \ref{mp5-inj.sym2}, there exist $M_1, M_2>0$ such that
 for all $v\in T$, $|\phi(v)|\leq |v|+ M_1$ and $|\phi^{-1}(v)|\leq |v|+ M_2$. Since $\phi$ is bijective, we have,
 $|\phi(v)|\leq |v|+ M_1$ and $|v|\leq |\phi(v)|+ M_2$ for all $v\in T$. By taking $M=\max\{M_1,M_2\}$, we get
 $|\,|\phi(v)|-|v|\,|\leq M$ for all $v\in T$.

 For the converse part, assume $\phi$
is bijective and that there exists an $M>0$ such that $|\,|\phi(v)|-|v|\,|\leq M$ for all $v\in T$.
This gives that $|\phi(v)|\leq |v|+ M$ and $|v|\leq |\phi(v)|+ M$ for all $v\in T$. Equivalently,
$|\phi(v)|\leq |v|+ M$ and $|\phi^{-1}(v)|\leq |v|+ M$ for all $v\in T$. Then,  by Corollary \ref{mp5-inj.sym2},
we get that both $C_\phi$, $C_{\phi^{-1}}$ are bounded operators on $\mathbb{T}_p$. Thus, $C_\phi$ is an invertible
 operator on $\mathbb{T}_p$ with an inverse $C_{\phi^{-1}}$.
\epf
\section{Isometry} \label{MP5Sec5}

Recall that a bounded operator $A$ on a normed linear space $(X,\|.\|)$ is said to be an \textit{isometry}
if $\|Ax\| = \|x\|$ for all $x\in X$.

\bthm
$C_\phi$ is an isometry on $\mathbb{T}_\infty$  if and only if $\phi$ is onto.
\ethm
\bpf
Suppose  that $\phi$ is onto. Then, since $\phi(T)=T$, $\|f\circ \phi\|_\infty = \|f\|_\infty$ for all $f\in \mathbb{T}_\infty$,
and hence $C_\phi$ is an isometry on $\mathbb{T}_\infty$.

Conversely, assume that $C_\phi$ is an isometry on $\mathbb{T}_\infty$. Now, suppose on the contrary that $\phi$ is not onto. Then, choose
$w\notin \phi(T)$. If $f=\chi_{w}$, then $\|f\circ \phi\|_\infty \neq \|f\|_\infty$, which is a contradiction. The result follows.
\epf
\bcor
$C_\phi$ is an isometry on $\mathbb{T}_{\infty,0}$  if and only if $\phi$ is onto and
$|\phi(v)|\rightarrow\infty$ as $|v|\rightarrow\infty$.
\ecor
\bthm\label{th-isometry1}
Let $T$ be a $2$-homogeneous tree and let $1\leq p< \infty$. Then  $C_\phi$ is an isometry on $\mathbb{T}_p$  if and only if
the following properties hold:
\begin{enumerate}
\item $\phi(\textsl{o})=\textsl{o}$
\item $\phi$ is onto
\item $|\phi(v)|=|\phi(w)|$ whenever $|v|=|w|$
\item If $\phi(w)\neq \textsl{o}$ for some $w\in T$, then $\phi$ is injective on $D_{|w|}$.
\end{enumerate}
\ethm
\bpf
Assume that $C_\phi$ is an isometry on $\mathbb{T}_p$.

First let us suppose that $\phi(\textsl{o})\neq\textsl{o}$.  If $f=\chi_{\textsl{o}}+(2)^{\frac{1}{p}}\chi_{\phi(\textsl{o})}$, then
$\|f\circ \phi\|\neq \|f\|$, which is a contradiction. Thus $(1)$ holds.

Secondly, let us suppose that $\phi$ is not onto. Then pick a $w\notin \phi(T)$. If $f=\chi_{w}$, then $\|f\circ \phi\|\neq \|f\|$,
which is again a contradiction. So, $(2)$ holds.

Thirdly, let us assume that there exist $v_1, v_2 \in T$ such that $|v_1|=|v_2|$
and $|\phi(v_1)|\neq |\phi(v_2)|$. Let $w_1=\phi(v_1)$
and  $w_2=\phi(v_2)$. Then
take
$$ f= (c_{|w_1|})^{\frac{1}{p}}\chi_{w_1}+(c_{|w_2|})^{\frac{1}{p}}\chi_{w_2},
$$
and observe that $\|f\|=1$. But,
$$\|C_\phi(f)\|^p \geq M_p^p(|v_1|,C_\phi(f))\geq 3/2,
$$
which is not possible. Thus property (3) holds.

Finally, let us suppose that there exists a $v_1\in T$ such that $\phi(v_1)\neq \textsl{o}$ and $\phi$ is not injective on $D_{|v_1|}$.
By property $(3)$, $\phi\not\equiv\textsl{o}$ on $D_{|v_1|}$. Since $\phi$ is not injective on $D_{|v_1|}$, we have
$\phi(v_1)=\phi(v_2)=w$ (say), where $|v_1|=|v_2|$. Now, we
take $f= 2^{\frac{1}{p}}\chi_{w}$. Then, $\|f\|=1$. But,
$$\|C_\phi(f)\|^p \geq M_p^p(|v_1|,C_\phi(f))=2,
$$
which is a contradiction, and hence $(4)$ holds.

Conversely, assume that all the four properties $(1)-(4)$ hold. We need to show that $C_\phi$ is an isometry on $\mathbb{T}_{p}$.
To do this, we fix $f\in \mathbb{T}_{p}$. Then, by the property (1), we have
$$|f(\phi(\textsl{o}))|^p=|f(\textsl{o})|^p \leq \|f\|^p .
$$
Fix $n\in\N$. By properties $(3)$ and $(4)$, $\phi \equiv 0$ on $D_n$ or $\phi$ is bijective from $D_n$ onto
$D_m$ for some $m\in \N$. In either case, $M_p^p(n,C_\phi(f))\leq \|f\|^p$, and thus
$$\|C_\phi(f)\|^p\leq \|f\|^p.
$$
Now, fix $m\in\N$. By properties $(2)$ and $(4)$, there exists an $n_m\in\N$ such that $\phi$ maps bijectively from $D_{n_m}$ onto
$D_m$. Therefore,
$$ M_p^p(m,f)= M_p^p(n_m,C_\phi(f))\leq \|f\circ \phi\|^p,
$$
and thus $\|f\|^p \leq \|C_\phi(f)\|^p $. Hence, $C_\phi$ is an isometry on $\mathbb{T}_{p}$.
\epf

\bcor
Let $1\leq p< \infty$.
The operator $C_\phi$ is an isometry on $\mathbb{T}_{p,0}$ over $2$-homogeneous trees if and only if
$\phi$ induces an isometric composition operator $C_\phi$ on $\mathbb{T}_{p}$, and $|\phi(v)|\rightarrow\infty$ as $|v|\rightarrow\infty$.

\ecor
\brem
If $C_\phi$ is an isometry on $\mathbb{T}_{p}$ (or $\mathbb{T}_{p,0}$) over $2$-homogeneous trees,
then the properties $(3)$ and $(4)$  in  Theorem \ref{th-isometry1} hold. However,
this is not the case for $(q+1)$-homogeneous trees with $q\geq2$. We will provide an example to demonstrate this fact.

For $v\neq \textsl{o}$, $v^-$ denote the parent of $v$. Fix $k\in \N$ with $k\geq 2$. For each element of $D_{k-1}$, choose
a child vertex and denote them by $v_1,v_2,\ldots,v_{c_{k-1}}$. Define
$$ \phi(v)= \left \{\begin{array}{ll}
v & \mbox{ if }~ |v|<k,\\
\textsl{o} &\mbox{ if }~ v\in D_k, v\neq v_1,v_2,\ldots,v_{c_{k-1}},\\
v^- & \mbox{ elsewhere.}
\end{array} \right.
$$
For $f\in \mathbb{T}_{p}$, we easily see that,
$$ M_p^p(n,C_\phi(f))= \left \{\begin{array}{ll}
M_p^p(n,f) & \mbox{ if }~ n<k,\\
M_p^p(n-1,f) &\mbox{ if }~ n>k,
\end{array} \right.
$$
and $ M_p^p(k,C_\phi(f)) \leq \|f\|^p$. This gives that, $C_\phi$ is an isometry on $\mathbb{T}_{p}$ (or $\mathbb{T}_{p,0}$).

Note that some vertices of $D_k$ are mapped into its parents, but all other vertices in $D_k$ are mapped to $\textsl{o}$. Consequently,
$\phi$ violates both the properties $(3)$ and $(4)$. The desired assertion follows.
\erem

It is natural to characterize the  isometric composition operators $C_\phi$ over $(q+1)$-homogeneous trees with $q\geq2$.

\bthm\label{th-isometry2}
Let $T$ be a $(q+1)$-homogeneous tree with $q\geq2$ and let $1\leq p< \infty$. Denote $\frac{c_k N_{k,n}}{c_n}$
by $\lambda_{k,n}$. Then, $C_\phi$ is an isometry on $\mathbb{T}_p$  if and only if the following properties hold:
\begin{enumerate}
\item $|\phi(v)|\leq |v|$. In particular, $\phi(\textsl{o})=\textsl{o}$.
\item $\sum\limits_{k=0}^n \lambda_{k,n}=1$ for all $n\in \N_0$.
\item For each $k\in \N_0$, $N_{\phi}(n,w)=N_{k,n}$ whenever $|w|=k$.
\item  $\sup\limits_{n\in \N_0}\lambda_{k,n}=1$ for all $k\in \N_0$. In particular, $\phi$ is onto.
\end{enumerate}
\ethm
\bpf
Assume that $C_\phi$ is an isometry on $\mathbb{T}_p$.

Suppose there exists a $v\in T$ such that $|v|<|\phi(v)|$. Let $w=\phi(v)$. Then the function
$f=(c_{|w|})^{\frac{1}{p}}\chi_{w}$ belonging to $\mathbb{T}_p$ contradicts the fact that $C_\phi$ is an isometry. Hence,
property $(1)$ holds.

Fix $n\in \N_0$. By property $(1)$, $\phi(D_n)\subseteq \bigcup\limits_{m=0}^n D_m$.
For each $k=0,1,2,\ldots,n$, choose $v_k\in D_k$ such that $N_{k,n}=N_\phi(n,v_k)$. Take
$f=\sum\limits_{k=0}^n(c_{k})^{\frac{1}{p}}\chi_{v_k}$ so that
$$\frac{1}{c_n}\sum\limits_{k=0}^n c_k N_{k,n}=M_p^p(n,f\circ\phi)\leq  \|f\circ\phi\|^p=\|f\|^p=1, ~\mbox{ i.e., }~\sum\limits_{k=0}^n c_k N_{k,n} \leq c_n.
$$
%which implies that
%$$\sum\limits_{k=0}^n c_k N_{k,n} \leq c_n.
%$$
By the definition of $N_{k,n}$, one can note that the number of vertices in $D_n$ which are mapped into $D_k$ under $\phi$ is less than or
equal to $c_k N_{k,n}$. Again by $(1)$, we get
$$c_n\leq \sum\limits_{k=0}^n c_k N_{k,n}.
$$
Therefore, $\sum\limits_{k=0}^n \lambda_{k,n}=1$ for all $n\in \N_0$. Thus property $(2)$ holds.

Suppose that there exist an $n_1\in\N$, and a $w\in T$ such that $N_{\phi}(n_1,w)<N_{k,n_1}$. Therefore, the
number of vertices in $D_{n_1}$ which are mapped into $D_k$ is strictly less than $c_k N_{k,n_1}$. Then by $(2)$, the total number of
elements in $D_{n_1}$ is strictly less than
$$\sum\limits_{m=0}^{n_1}c_m N_{m,n_1}=c_{n_1},
$$
which is a contradiction. Thus, the property $(3)$ is verified.

Fix $k\in \N_0$ and $w\in D_k$. Take $f=(c_{|w|})^{\frac{1}{p}}\chi_{w}$. By $(3)$, for each $n\in \N_0$, we see that
$$M_p^p(n,f\circ\phi)=\frac{c_{|w|}}{c_n}N_\phi(n,w)=\frac{c_k N_{k,n}}{c_n}=\lambda_{k,n}
$$
and hence,
$$\sup\limits_{n\in \N_0}\lambda_{k,n}=\|f\circ\phi\|^p=\|f\|^p=1.
$$

Conversely, assume that all the four properties $(1)-(4)$ hold.
In order to prove that  $C_\phi$ is an isometry on $\mathbb{T}_p$, we fix $f\in \mathbb{T}_{p}$. Then,  for $n\in \N_0$, we have
\beqq
M_p^p(n,C_\phi f)
&=& \frac{1}{c_n}  \sum\limits_{k=0}^{n} \sum_{\substack{|\phi(v)|=k \\ |v|=n}} |f(\phi (v))|^p ~\mbox{ (by (1))}\\
&=&  \frac{1}{c_n}\sum\limits_{k=0}^{n} c_k N_{k,n} M_p^p(k,f)  ~\mbox{ (by (3))}\\
&=& \sum\limits_{k=0}^{n} \lambda_{k,n} M_p^p(k,f).
\eeqq
Thus, $\|C_\phi f\|^p= \sup A$, where
$$A:=\left\{s_n=\sum\limits_{k=0}^{n} \lambda_{k,n} M_p^p(k,f): n\in \N_0\right\}.
$$
For each $n\in \N_0$, $s_n\leq \|f\|^p$ which implies $\sup A\leq\|f\|^p$.
Fix $m\in \N_0$. Then, by $(4)$, we have $\sup\limits_{n\in \N_0}\lambda_{m,n}=1$ and that
there exists  an $n_1\in \N_0$ such that $\lambda_{m,n_1}=1$, or else, there exists a subsequence
$\{\lambda_{m,n_k}\}$ converging to $1$ as $k\rightarrow\infty$. In the first case, $s_{n_1}=M_p^p(m,f)\in A$
so that $M_p^p(m,f)\leq \sup A$. In the later case,
\beqq
|M_p^p(m,f)-s_{n_k}|
&\leq& \sum_{\substack{k=0 \\ k\neq m}}^{n_k} \lambda_{k,n_k} M_p^p(k,f)+ (1-\lambda_{m,n_k}) M_p^p(m,f)\\
&\leq& 2(1-\lambda_{m,n_k}) \|f\|^p,
\eeqq
which implies that $M_p^p(m,f)$ is a limit point of $A$ and therefore, $M_p^p(m,f)\leq \sup A$. Since $m$ was
arbitrary, we have $\|f\|^p\leq \sup A$. Hence, $\|f\|^p=\sup A=\|f\circ\phi\|^p$. The desired conclusion follows.
\epf

\bcor
Let $1\leq p<\infty$  and $C_\phi$ be a bounded operator on $\mathbb{T}_{p,0}$ over $(q+1)$-homogeneous tree with
$q\geq2$. Then $C_\phi$ is an isometry on  $\mathbb{T}_{p,0}$ if and only if $C_\phi$ is an isometry on  $\mathbb{T}_{p}$.
\ecor
\bcor
Let  $1\leq p< \infty$. Suppose $C_\phi$ is an isometry on $\mathbb{T}_p$ or on $\mathbb{T}_{p,0}$,  and $|\phi(v)|=|v|$ for some
$v\neq \textsl{o}$. Then $\phi$ is a permutation on $D_{|v|}$.
\ecor
\bpf
The results holds easily for all $2$-homogeneous trees and thus, we assume that $T$ to be a $(q+1)$-homogeneous
tree with $q\geq2$. Assume that there exist  $v_1, v_2 \in T$ with $|v_1|=|v_2|\neq0, |\phi(v_1)|=|v_1|$
and $|\phi(v_2)|\neq|v_2|$. Consider the function
$$f= (c_{|v_1|})^{\frac{1}{p}}\chi_{w_1}+(c_{|w_2|})^{\frac{1}{p}}\chi_{w_2},
$$
where  $w_1=\phi(v_1)$  and $w_2=\phi(v_2)$.
Since $|w_1|\neq|w_2|$, we have $\|f\|=1$. But,
$$\|C_\phi(f)\|^p \geq M_p^p(|v_1|,C_\phi(f))> 1,
$$
which contradicts the hypothesis. Thus, $\phi(D_{|v_1|})\subseteq D_{|v_1|}$.

Next, we claim that $\phi$ is a permutation on $D_{|v_1|}$. Suppose not. Then there exist $w_1, w_2 \in D_{|v_1|}$ with
$\phi(w_1)=\phi(w_2)=w$ (say). Thus $C_\phi$ is not an isometry via the test function $f= (c_{|v_1|})^{\frac{1}{p}}\chi_{w}$.
Thus, $\phi$ is injective on $D_{|v_1|}$ and the result follows.
\epf

\bcor
Suppose $\phi \in \Aut (T)$. Then, $C_\phi$ is an isometry on $\mathbb{T}_{p}$ if and only if
$\phi(\textsl{o})= \textsl{o}$.
\ecor
\bpf
Suppose that $\phi \in \Aut (T)$ and $\phi(\textsl{o})= \textsl{o}$. Then, $\phi$ is a bijective
map from $D_n$ to $D_n$ for each $n$. Therefore, for $n\in \mathbb{N}_0$ and
$f\in \mathbb{T}_{p}$, we  have $M_p^p(n,C_\phi f) = M_p^p(n,f)$. Hence,
$C_\phi$ is an isometry on $\mathbb{T}_{p}$.

Converse part is already contained in Theorems \ref{th-isometry1} and \ref{th-isometry2}.
\epf

\section{Compact Composition Operators} \label{MP5Sec6}

Recall that a bounded operator $A$ on a normed linear space $X$ is said to be
\textit{compact operator} if the image of the closed unit ball is relatively compact in $X$.
We now recall some important results from \cite{CO-Tp-spaces}.

\begin{Thm}%\label{}
{\rm (\cite[Corollary 5]{CO-Tp-spaces})}
Let $\phi$ be a self-map of $T$. Then $C_\phi$ is compact on $\mathbb{T}_{\infty}$
if and only if $\phi$ is a bounded self-map of $T$.
\end{Thm}

\begin{Thm}%\label{MP2-th7}
{\rm (\cite[Corollary 8]{CO-Tp-spaces})}
Let $T$ be a $2$-homogeneous tree. Then $C_\phi$ is compact on $\mathbb{T}_{p}$ if and only if
$\phi$ is a bounded self-map of $T$.
\end{Thm}

\begin{Thm}\label{MP2-th6}
{\rm (\cite[Theorem 7]{CO-Tp-spaces})}
If  $\phi$ is a self-map of $(q+1)$-homogeneous tree $T$ with $q\geq 2$, then the following are equivalent:
\begin{enumerate}
\item[{\rm (a)}] $C_\phi$ is compact on $\mathbb{T}_{p}$ for $1\leq p\leq\infty$.
\item[{\rm (b)}] $\|C_\phi f_n\| \rightarrow 0$ as $n\rightarrow \infty$ whenever $\{f_n\}$ is a bounded sequence of functions
that converges to $0$ pointwise.
\end{enumerate}
\end{Thm}

\begin{Thm}%\label{MP2-cor6}
{\rm (\cite[Corollary 7]{CO-Tp-spaces})}
If $C_\phi$ is compact on $\mathbb{T}_{p}$, then $|v|-|\phi(v)|\rightarrow \infty$ as $|v|\rightarrow\infty$.
\end{Thm}

\bthm%\label{norm}
Let $T$ be a  $(q+1)$-homogeneous tree with $q\geq 2$, and   $1\leq p<\infty$. Then $C_\phi$ is compact operator
on $\mathbb{T}_{p}$ whenever
$$\frac{1}{c_n}\sum\limits_{m=0}^{\infty}N_{m,n}c_m \rightarrow 0 ~\mbox{ as }~ n\rightarrow \infty.
$$
\ethm
\bpf
Let $\{f_k\}$ be a bounded sequence such that $\{f_k\}$  converges to $0$ pointwise. Without loss of
generality, we may assume that $\|f_k\|\leq 1$ for all $k$. Then, by Theorem \Ref{MP2-th6}, it suffices to show
that $\|C_\phi (f_k)\|\rightarrow 0$ as $k\rightarrow\infty$.

Fix $\epsilon>0$. Then, by the hypothesis, there exists an $N_1\in \N$ such that
\be\label{MP5-eq2}
\frac{1}{c_n}\sum\limits_{m=0}^{\infty}N_{m,n}c_m \leq \epsilon^p ~\mbox{ for all }~ n\geq N_1.
\ee
Set $S=\{\phi(v): |v|<N_1\}$. Then,  since $\{f_k\}$ converges to $0$ pointwise and $S$ is a finite set, it follows that
$\{f_k\}$ converges to $0$ uniformly on $S$ and thus, there exists an $N\in \N$ such that
$$\sup \limits_{w\in S} |f_k(w)|\leq \epsilon ~\mbox{ for all }~ k\geq N.
$$
Fix $k\geq N$. Then, for $n\in \N_0$ with $n<N_1$, we have $M_p^p(n,C_\phi f_k) \leq \epsilon^p$. Next, for $n\geq N_1$, we have
\beqq
M_p^p(n,C_\phi f_k)
&=& \frac{1}{c_n}  \sum\limits_{m=0}^{\infty} \sum_{\substack{|\phi(v)|=m \\ |v|=n}} |f_k(\phi (v))|^p \\
&\leq& \frac{1}{c_n}\sum\limits_{m=0}^{\infty} c_m N_{m,n} \|f_k\|^p  \\
&\leq&  \epsilon^p ~\mbox{ (by \eqref{MP5-eq2})}
\eeqq
which shows that $\|C_\phi f_k\|\rightarrow 0$ as $k\rightarrow\infty$. Thus, $C_\phi$ is compact on
$\mathbb{T}_{p}$.%, by Theorem \ref{MP2-th6}.
\epf

\bthm\label{Tp0-cmpt}
If $C_\phi$ is compact on $\mathbb{T}_{p,0}$ for $1\leq p\leq\infty$, then
$\|C_\phi f_n\| \rightarrow 0$ as $n\rightarrow \infty$ whenever $\{\|f_n\|: n\in \N\}$ is bounded,
and $\{f_n\}$ converges to $0$ pointwise.
\ethm
\bpf
Proof of this result is similar to the proof of the implication ``(a)~$\Ra$~(b)" in Theorem \Ref{MP2-th6}.
So we omit the details.
\epf

\bthm
There are no compact composition operators on  $\mathbb{T}_{\infty,0}$.
\ethm
\bpf
Suppose that $\phi$ is a bounded self-map of $T$. Then, by Theorem \ref{bdd:T_infty,0}, $C_\phi$ is not
bounded and hence it is not compact. Suppose that $\phi$ is an unbounded self-map of $T$.
Then, there exists a sequence of vertices $\{v_n\}$  such that
$\phi(v_n)=w_n$ and $|w_n|\rightarrow\infty$ as $n\rightarrow\infty$. Take $f_n=\chi_{\{w_n\}}$ for each $n\in\N$.
Then it easy to see that $\{f_n\}$ converges to $0$ pointwise and $\|C_\phi(f_n)\|_\infty=\|f_n\|_\infty=1$
for each $n$. Therefore, $C_\phi$ cannot be a compact operator on $\mathbb{T}_{\infty,0}$ by Theorem \ref{Tp0-cmpt}.
\epf
\bthm
Let $T$ be a $2$-homogeneous tree and $1\leq p<\infty$. Then, there are no compact composition
operators on  $\mathbb{T}_{p,0}$.
\ethm
\bpf
By Theorem \ref{bdd:Tp,0,q=1}, every bounded self-map of $T$ cannot induce a bounded (in particular, compact)
composition operator. Suppose that $\phi$ is an unbounded self-map of $T$. Then, choose a sequence of vertices $\{v_n\}$
such that $\{w_n\}$ is unbounded, where $\phi(v_n)=w_n$.
Take $f_n=2^{1/p}\chi_{\{w_n\}}$ so that $\{f_n\}$ converges to $0$ pointwise and $\|f_n\|=1$
for each $n$. Finally, since
$$1\leq M_p^p(|v_n|,f_n\circ\phi)\leq \|C_\phi(f_n)\|^p ~\mbox{ for all }~ n\in \N,
$$
it follows that $C_\phi$ cannot be a compact operator on $\mathbb{T}_{p,0}$ by Theorem \ref{Tp0-cmpt}.
\epf
\bthm
Let $T$ be a $(q+1)$-homogeneous tree with $q\geq 2$. Then the operator $C_\phi$ cannot be compact on
 $\mathbb{T}_{p,0}$, $1\leq p<\infty$, for any self-map $\phi$ of $T$.
\ethm
\bpf
 Suppose $C_\phi$ is compact for a self-map $\phi$ of $T$. Consider the sequence of functions
 defined by
 $$g_n(v)=\frac{n}{n+|v|} ~\mbox{ for }~ v\in T, n\in\N.
 $$ It is easy to see that
 $$
 M_p(m,g_n)=\frac{n}{n+m} ~\mbox{ for }~ n\in \N, m\in \N_0.
 $$
Therefore, $g_n\in \mathbb{T}_{p,0}$ with $\|g_n\|=1$ for all $n\in\N$. For each fixed $v\in T$,
$g_n(v)\ra 1$ pointwise. Since $C_\phi$ is compact on $\mathbb{T}_{p,0}$, there exists a subsequence
$\{g_{n_k}\}$ of $\{g_{n}\}$ and $g\in \mathbb{T}_{p,0}$ such that $C_\phi(g_{n_k})\ra g$ in $\|.\|_p$.
Then, by Lemma \Ref{lem:bound}, we have $g_{n_k}(\phi(v))\ra g(v)$ pointwise for $v\in T$,
which gives that $g\equiv1$.
Since $g\notin \mathbb{T}_{p,0}$, $C_\phi$ cannot be compact on $\mathbb{T}_{p,0}$.
\epf

%\subsection*{\bf Acknowledgement}
%%The authors thank the referee for many useful comments.
%The first author thanks...% the Council of Scientific and Industrial Research (CSIR), India,
%%for providing financial support in the form of a SPM Fellowship to carry out this research.
%%The second  author is currently on leave from the IIT Madras.
%%

\subsection*{Acknowledgments}
The authors thank the referee for many useful comments which improve the
presentation considerably.
%The  work of the second author is supported by Mathematical Research Impact Centric Support (MATRICS) of DST, India  (MTR/2017/000367).
%,and the second author is currently at ISI, Chennai Centre, Chennai, India.
%The authors thank the referee for his/her comments.

\subsection*{Conflict of Interests}
The authors declare that there is no conflict of interests regarding the publication of this paper.

\end{document}